\documentclass[aos,MSNbibl,nameyear,dvips]{arximspdf}
\usepackage{dcolumn}
\usepackage{graphicx}

\doi{10.1214/12-AOS982} 
\volume{40}
\issue{2}
\pubyear{2012}
\firstpage{812}
\lastpage{831}

\makeatletter
\newcolumntype{d}[1]{D{.}{.}{#1}}
\newtheorem{lemma}{Lemma}
\newtheorem{theorem}{Theorem}
\def\back{}
\def\bfX{\mathbf{X}}
\def\bfl{\mathbf{l}}
\def\bfA{\mathbf{A}}
\def\bfB{\mathbf{B}}
\def\bfy{\mathbf{y}}
\def\bfD{\mathbf{D}}
\def\bfP{\mathbf{P}}
\def\bfQ{\mathbf{Q}}
\def\bfI{\mathbf{I}}
\def\bftheta{\bolds{\theta}}
\def\bfvartheta{\bolds{\vartheta}}
\def\hatbftheta{\hat{\bolds{\theta}}}
\def\hatbfvartheta{\hat{\bolds{\vartheta}}}
\def\bfGamma{\bolds{\Gamma}}
\def\bfbeta{\bolds{\beta}}
\def\bfxi{\bolds{\xi}}
\def\bfzeta{\bolds{\zeta}}
\def\bfSigma{\bolds{\Sigma}}
\def\bfzero{\mathbf{0}}
\def\calM{\mathcal{M}}
\makeatother

\begin{document}
\begin{frontmatter}

\title{Estimation in high-dimensional linear models with~deterministic design matrices}
\runtitle{Estimation in high-dimensional linear models}

\begin{aug}
\author[A]{\fnms{Jun} \snm{Shao}\corref{}\ead[label=e1]{shao@stat.wisc.edu}\thanksref{t1}}
\and
\author[B]{\fnms{Xinwei} \snm{Deng}\ead[label=e2]{xdeng@vt.edu}}
\thankstext{t1}{Supported in part by NSF Grants SES-07-05033 and DMS-10-07454.}
\runauthor{J. Shao and X. Deng}
\affiliation{East China Normal University, University of Wisconsin
and
Virginia~Polytechnic Institute and State University}
\address[A]{School of Finance and Statistics\\
East China Normal University\\
500 Dongchuan Rd.\\
Shanghai, 200241\\
China\\
and\\
Department of Statistics\\
University of Wisconsin\\
1300 University Ave.\\
Madison, Wisconsin 53706\\
USA\\
\printead{e1}} 
\address[B]{Department of Statistics\\
Virginia Polytechnic Institute\\
\quad and State University\\
211 Hutcheson Hall\\
Blacksburg, Virginia 24061\\
USA\\
\printead{e2}}
\end{aug}

\received{\smonth{9} \syear{2010}}
\revised{\smonth{1} \syear{2012}}

%
\begin{abstract}
Because of the advance in technologies,
modern statistical studies often encounter linear models
with the number of explanatory variables much larger than
the sample size.
Estimation and variable selection in these high-dimensional
problems with deterministic design points
is very different from those in the case
of random covariates, due to the
identifiability of the high-dimensional
regression parameter vector. We show that a reasonable
approach is to focus on the projection of the regression parameter vector
onto the linear space generated by the design matrix.
In this work, we consider the ridge regression estimator of the
projection vector
and propose to threshold the ridge regression estimator when the
projection vector is sparse in the sense that many of its components
are small. The proposed estimator has an explicit form and
is easy to use in application.
Asymptotic properties such as the consistency of
variable selection and estimation and the
convergence rate of the prediction mean squared error are established
under some sparsity conditions on the projection
vector. A simulation study is also conducted to examine the
performance of the proposed estimator.
\end{abstract}

%


\end{frontmatter}

\section{Introduction}\label{sec1}

Consider the following linear model:
%
%
\begin{equation}
y_i = \mathbf{x}_i' \bfbeta+ \varepsilon_i ,\qquad  i=1,\ldots,n,
\label{model}
\end{equation}
where $y_i$ is an observed response variable,
$\mathbf{x}_i$ is a $p$-dimensional vector of observed covariates or
design points associated with $y_i$,
$\bfbeta$ is a $p$-dimensional vector of unknown
parameters and $\varepsilon_i$'s are independent and
identically distributed unobserved random errors
with mean 0 and unknown variance $\sigma^2$.
The theory of linear models is well established
for traditional applications where the dimension\vadjust{\goodbreak}
$p$ is fixed and the sample size $n >p$.
With modern technologies, however,
in many biological, medical, social and economical studies,
$p$ is comparable with, or much larger than, $n$, and making valid
statistical inference
is a great challenge.

In the case of $p < n$,
there is a rich literature on variable selection,
that is, identifying nonzero components of~$\bfbeta$ in (\ref{model}).
For variable selection in the case of $p>n$ and
statistical inference afterwards,
the development of statistical theory started about a decade ago.
Some excellent advances in asymptotic theory
have been made recently in situations where
$p$ diverges to
infinity as the sample size $n$ increases to infinity with
the divergence rate $O(n^l )$ for some $l > 0$
(polynomial-type divergence rate) or
$O(e^{n^\nu} )$ for some $\nu\in(0,1)$
(ultra-high dimension).
See, for example,
Fan and Peng (\citeyear{FanPen04}), Hunter and Li
(\citeyear{HunLi05}), Meinshausen and Buhlmann (\citeyear{MeiBuh06}),
Zhao and Yu (\citeyear{ZhaYu06}), Zou (\citeyear{Zou06}), Wang, Li and Tsai (\citeyear{WanLiTsa07}),
Fan and Lv (\citeyear{FanLv08}), Zhang and Huang (\citeyear{ZhaHua08}),
Meinshausen and Yu (\citeyear{MeiYu09}), Wang (\citeyear{Wan09})
and a review by Fan and Lv (\citeyear{FanLv10}).
When $\mathbf{x}_i$'s are random covariates,
under some conditions, some variable selection methods have
been shown to be selection-consistent in the sense that,
with probability tending to 1 as $n \rightarrow
\infty$, the selected variables are exactly those
related to the response, where the probability is
with respect to the joint distribution of $(y_i,\mathbf{x}_i)$'s.
As Fan and Lv (\citeyear{FanLv08}) commented in the end of
their stimulating paper, however, no selection-consistency result
is available for deterministic $\mathbf{x}_i$'s and many applications,
such as
biomedical imaging and signal processing, involve deterministic design points.
Another example in which $\mathbf{x}_i$ can be treated as deterministic is
an analysis conditional on the observed covariates.

Let $\bfX$ be the matrix whose $i$th row is $\mathbf{x}_i'$, $i=1,\ldots,n$.
For simplicity, we call $\bfX$ the design matrix although $\mathbf{x}_i$'s
are not necessarily designed points. When $p>n$, a~key difference
between a random $\bfX$ and a deterministic design matrix is the
identifiability of the regression parameter~$\bfbeta$ in
(\ref{model}), caused by the fact that the probabilities under
consideration are different. For random $\mathbf{x}_i$'s that are
independent and identically distributed and independent of
$\varepsilon_i$'s, $\bfbeta= [\operatorname{cov} (\mathbf
{x}_i)]^{-1}\operatorname{cov}
(\mathbf{x}_i , y_i)$. Hence, even when $p>n$, components of~$\bfbeta$ can
be estimated, and nonzero components of~$\bfbeta$ can be identified
consistently with respect to the joint probability distribution of
$(y_i,\mathbf{x}_i)$'s, under some conditions on $\operatorname{cov} (\mathbf
{x}_i)$ and
$\operatorname{cov} (\mathbf{x}_i , y_i)$. On the other hand, when the design
matrix is deterministic or an analysis conditional on $\bfX$ is
considered, the underlying probability is the probability
distribution of $(y_1,\ldots,y_n)$ conditional on $\bfX$, and
$\bfbeta$
is identifiable if and only if it lies in a set having a~one-to-one
correspondence with $\mathcal{R} ( \bfX)$, the linear space spanned by
rows of~$\bfX$. Since the dimension of $\mathcal{R}(\bfX)$ is at most
$n$, when $p>n$,~$\bfbeta$ is generally not identifiable with
respect to the probability distribution of $(y_1,\ldots,y_n)$
conditional on $\bfX$. Consequently, with deterministic $\bfX$ and
$p >n$, it is not realistic to derive consistent estimators of
$\bfbeta$ or consistent variable selection procedures.

Without selection-consistency [as previously described; see
definition (\ref{sc}) in Section~\ref{sec4.1}],
we may still\vadjust{\goodbreak}
derive consistent estimators of some useful functions of~$\bfbeta$
under the $p$-dimensional linear model
given by (\ref{model}) with deterministic $\bfX$ and $p>n$.
This is the main focus of the current paper.
Although~$\bfbeta$ is generally not
identifiable when $p>n$, we argue in Section~\ref{sec2} that
we may not need to estimate the entire vector~$\bfbeta$.
For statistical analysis, $\bftheta$,
the projection of~$\bfbeta$ onto $\mathcal{R} ( \bfX)$,
is what we are able to estimate, and perhaps the estimation
of $\bftheta$ is sufficient for valid statistical inference.

To estimate $\bftheta$, we first consider the ridge
regression estimator in Section~\ref{sec3}.
For any linear combination of the ridge regression estimator,
we establish the asymptotic convergence rate of its mean squared error.
We also obtain the convergence rate of
the expected $L_2$-norm error for the ridge regression estimator of~$\bfX\bftheta$.
This expected $L_2$-norm error divided by $n$ is equal to the average
prediction mean squared error minus $\sigma^2$.

When $\bftheta$ is sparse in the sense that many of its components are small,
we consider in Section~\ref{sec4} a sparse estimator of $\bftheta$ obtained by
thresholding
the ridge regression estimator of $\bftheta$.
We show that, with probability tending to 1 at a fast rate,
we can eliminate small components of $\bftheta$ and keep large components
of $\bftheta$, that is, thresholding
the ridge regression estimator provides
a variable selection procedure, that is, consistent in some
sense. This method is computationally much simpler
than methods such as the LASSO [Tibshirani (\citeyear{Tib96})],
SCAD [Fan and Li (\citeyear{FanLi01})]
and the ENET [Zou and Hastie (\citeyear{ZouHas05})],
since no numerical minimization is required
as the proposed estimator has an explicit form. We show that the
convergence rate of
the expected $L_2$-norm error or average prediction mean squared error
of the thresholded ridge regression estimator is much faster than
that of the ridge regression estimator when $\bftheta$ is sparse.
In particular, the thresholded ridge regression estimator
is estimation-consistent (defined in Section~\ref{sec4}), but the
ridge regression estimator may not be.

Thresholding the ridge regression estimator
is closely related to the SIS as shown in Fan and Lv (\citeyear{FanLv08}).
However, its asymptotic behavior for deterministic~$\bfX$ is different from that for random $\bfX$, and its
consistency also requires very different conditions.
For deterministic $\bfX$ and $p>n$,
there does not exist any result on the consistency of the LASSO, SCAD
or ENET.
When $p < n$, Zhang and Huang (\citeyear{ZhaHua08}) showed that the LASSO is
estimation-consistent,
but the required conditions are more stringent and complicated
than those required for the consistency of
the thresholded ridge regression estimator.

Some simulation results are presented in Section~\ref{sec5}
to study the estimation and
prediction performance of the
proposed method, the ridge regression, the LASSO and
the ENET.
All technical proofs are given in Section~\ref{sec6}.

\section{Identifiability and projection}\label{sec2}

We consider model (\ref{model}) with deterministic
design matrix $\bfX= (\mathbf{x}_1,\ldots,\mathbf{x}_n)'$, where
the dimension of $\mathbf{x}_i$, $p $, is larger than $n$.
Let $r = r_n$ be the rank of $\bfX$.
From the singular value decomposition,
%
\begin{equation}
\bfX= \bfP\bfD\bfQ' , \label{decom}\vadjust{\goodbreak}
\end{equation}
where $\bfP$ is an $n \times r$ matrix satisfying $\bfP' \bfP= \bfI_r$,
$\bfQ$ is a $p\times r$ matrix satisfying $\bfQ' \bfQ= \bfI_r$,
$\bfI_a$ denotes the identity matrix of order $a$
and $\bfD$ is an $r \times r$ diagonal matrix of full rank.
Let $\bfQ_\perp$ be a $p \times(p-r)$ matrix such that
$\bfQ' \bfQ_\perp= \bfzero$
(the matrix of 0's with an appropriate order)
and $\bfQ_\perp' \bfQ_\perp= \bfI_{p-r}$.
Throughout, we denote the $q$-dimensional Euclidean space by $\mathcal{R}^q$
for any positive integer~$q$ and
the subspace of $\mathcal{R}^p$ generated by the rows of $\bfX$ by
$\mathcal{R} ( \bfX)$.

We say that~$\bfbeta$ in (\ref{model}) is identifiable if
$\bfbeta_1 \in\bfB$, $\bfbeta_2 \in\bfB$ and
$\bfX\bfbeta_1 = \bfX\bfbeta_2$ imply $\bfbeta_1 = \bfbeta_2$, where
$\bfB$ is the parameter space of~$\bfbeta$. The following lemma
gives a sufficient and necessary condition for the identifiability
of~$\bfbeta$.

\begin{lemma}\label{le1}
Under model (\ref{model}) with $p>r$,
$\bfbeta$ is identifiable if and only if there exists a known function
$\phi$
from $\mathcal{R}^r$ to $\mathcal{R}^{p-r}$ such that
%
%
\begin{equation}
\bfB= \{ \bfbeta\dvtx  \bfbeta= \bfQ\bfxi+ \bfQ_\perp\phi( \bfxi) ,
 \bfxi\in\mathcal{R}^r \}. \label{bspace}
\end{equation}
\end{lemma}

Lemma~\ref{le1} reveals that identifiable~$\bfbeta$'s
must be in a set having a one-to-one correspondence
with $\mathcal{R}( \bfX) = \{ \bfbeta\dvtx  \bfbeta= \bfQ\bfxi,
 \bfxi\in\mathcal{R}^r\}$. Since the dimension of the set on the
right-hand side
of (\ref{bspace}) is $r \leq n \wedge p $ (the minimum of~$n$ and~$p$),
$\bfbeta$ is typically not identifiable when $p> n $ and, hence, we are
not able to obtain
a component-wise consistent estimator of~$\bfbeta$.
However, we may not need to estimate the entire vector~$\bfbeta$, that is, if $\bfX\bfbeta_1 = \bfX\bfbeta_2$, we can
still estimate
parameters related to $\bfX\bfbeta_1 = \bfX\bfbeta_2$ and make valid
inference without trying to distinguish $\bfbeta_1$ and $\bfbeta_2$.
Therefore, we consider the projection
of~$\bfbeta$ onto~$\mathcal{R}( \bfX)$, which
is what we are able to identify in view of Lemma~\ref{le1}. Define
\[
(\bfX\bfX' )^- = \bfP\bfD^{-2} \bfP' ,
\]
which is $(\bfX\bfX' )^{-1}$ if $r=n$.
The projection of~$\bfbeta$ onto $\mathcal{R}(\bfX)$ is
%
%
\begin{equation}
\bftheta= \bfX' (\bfX\bfX' )^- \bfX\bfbeta= \bfQ\bfQ' \bfbeta
.\label{proj}
\end{equation}
Note that $\bftheta\in\mathcal{R}(\bfX)$ and
$\bftheta= \bfbeta$ if and only if $\bfbeta\in\mathcal{R}(\bfX)$.
Furthermore, $\bfX\bftheta= \bfX\bfbeta$ and model
(\ref{model}) can be written as
%
%
\begin{equation}
y_i = \mathbf{x}_i' \bftheta+ \varepsilon_i ,\qquad  i=1,\ldots,n.
\label{model1}
\end{equation}
Thus, estimating $\bftheta$ is enough for
inference about parameters
$\bfX\bfbeta= \bfX\bftheta$ and prediction.

The dimension of $\bftheta$ is still $p$.
When~$\bfbeta$ has many zero components, $\bftheta$ may not have any
zero component. However, $\bftheta$ may have many small components.
This can be seen from the $L_2$-norms of~$\bfbeta$ and $\bftheta$.
Since $\bftheta= \bfQ\bfQ' \bfbeta$ and
$\bfQ\bfQ'$ is a projection matrix, we obtain that
$ \| \bftheta\| \leq\| \bfbeta\| $, where $\| \cdot\|$ denotes the
$L_2$-norm.
This implies that if~$\bfbeta$ has many zero components so that the order
of $\| \bfbeta\|$ is
much smaller than $O(p)$, then the order of $\| \bftheta\|$ is
also much smaller than $O(p)$. Hence, if components of $\bftheta$ are nonzero,
then many of them must be negligible,
and $\bftheta$ can be viewed as a sparse vector. More precise descriptions
of this sparsity can be found in conditions~(\ref{eqc2}) in Section~\ref{sec3} and (\ref{eqc4})
in Section~\ref{sec4}.

\section{The ridge regression estimator of the projection}\label{sec3}

Since the dimension of $\bftheta$ in (\ref{proj}) is $p>n$, we consider
the ridge regression estimator of $\bftheta$
[Hoerl and Kennard (\citeyear{HoeKen70})] under model (\ref{model1}).
\[
\hatbftheta= ( \bfX' \bfX+ h_n \bfI_p )^{-1} \bfX' \bfy,
\]
where $\bfy= (y_1,\ldots,y_n)'$
and $h_n >0$ is an appropriately chosen regularization parameter.
The computation of $\hatbftheta$ involves only inverting
an $n \times n$ matrix. This is because (\ref{decom}) implies that
%
%
\begin{equation}
( \bfX' \bfX+ h_n \bfI_p )^{-1} \bfX'
= \bfX' ( \bfX\bfX' + h_n \bfI_n )^{-1} , \label{iden}
\end{equation}
which also implies that
the ridge regression estimator $\hatbftheta$
is always in $\mathcal{R}(\bfX)$.
In fact, if $ \hat{\bfbeta}$ is the ridge regression estimator of
$\bfbeta$
constructed under model~(\ref{model}), then $ \hatbftheta=
\bfX' (\bfX\bfX' )^- \bfX\hat{\bfbeta}= \hat{\bfbeta}$.
But $ \hatbftheta= \hat{\bfbeta}$ estimates
$\bftheta$, not the nonidentifiable~$\bfbeta$ when $p>n$.

We now study the bias and variance of $\hat{\bftheta}$
as an estimator of $\bftheta$, which is essential for
establishing asymptotic properties of $\hat{\bftheta}$.
For the matrix $\bfQ$ given in the singular value decomposition (\ref{decom}),
$\bfGamma= (  \bfQ  \bfQ_\perp )$ is
orthogonal, that is, $\bfGamma' \bfGamma= \bfGamma\bfGamma' = \bfI
_p$. Then
\begin{eqnarray*}
\operatorname{bias}( \hatbftheta) \back& = & \back E( \hat{\bftheta
} )
- \bftheta\\
& = & \back
( \bfX' \bfX+ h_n \bfI_p )^{-1} \bfX' \bfX\bftheta- \bftheta\\
& = & \back- ( h_n^{-1}\bfX' \bfX+ \bfI_p )^{-1} \bftheta\\
& = & \back- \bfGamma( h_n^{-1} \bfGamma' \bfX' \bfX\bfGamma
+ \bfI_p )^{-1} \bfGamma' \bfQ\bfQ' \bftheta\\
& = & \back- \pmatrix{\bfQ& \bfQ_\perp}
\pmatrix{
(h_n^{-1} \bfD^2 + \bfI_r )^{-1} & \bfzero\vspace*{2pt}\cr
\bfzero& \bfI_{p-r}}
\pmatrix{
\bfQ' \vspace*{2pt}\cr
\bfQ_\perp'}\bfQ\bfQ' \bftheta\\
& = & \back- \pmatrix{
\bfQ( h_n^{-1} \bfD^2 + \bfI_r )^{-1} & \bfQ_\perp
}
\pmatrix{\bfQ' \bftheta\vspace*{2pt}\cr
\bfzero}
\\
& = & \back- \bfQ
( h_n^{-1} \bfD^2 + \bfI_r )^{-1} \bfQ' \bftheta,
\end{eqnarray*}
where the fourth equality follows from the fact that
$\bfGamma$ is orthogonal and $\bftheta= \bfQ\bfQ' \bfbeta= \bfQ
\bfQ'
\bftheta$.
The covariance matrix of $\hat{\bftheta}$ is given by
\begin{eqnarray*}
\operatorname{var} ( \hatbftheta) \back& = & \back
\sigma^2 ( \bfX' \bfX+ h_n \bfI_p )^{-1}
\bfX' \bfX( \bfX' \bfX+ h_n \bfI_p )^{-1}\\
& \leq& \back\sigma^2 ( \bfX' \bfX+ h_n \bfI_p )^{-1} \\
& \leq& \back\sigma^2 h_n^{-1} \bfI_p ,
\end{eqnarray*}
where $\bfA\leq\bfB$ for nonnegative definite matrices $\bfA$ and
$\bfB$
means $\bfB- \bfA$ is nonnegative definite.

To study the asymptotic properties of
$\hatbftheta$, we consider $n \rightarrow\infty$ and $p=p_n$, a~function of $n$.
Quantities such as~$\bfbeta$, $\bfy$, $\mathbf{x}_i$, etc., form
triangular arrays, but the subscript $n$ is omitted for simplicity.
We assume that $\lambda_{1n}$, the
smallest positive eigenvalue of $\bfX' \bfX$, satisfies
\renewcommand{\theequation}{C\arabic{equation}}
\setcounter{equation}{0}
\begin{equation}\label{eqc1}
\lambda_{1n}^{-1} = O(n^{-\eta} ), \qquad
\eta\leq1 \mbox{ and }\eta\mbox{ does not depend on }n.
\end{equation}
We also need a sparsity condition on $\bftheta$. From the discussion
in the end of Section~\ref{sec2}, we conclude that, in terms of the $L_2$-norm,
the sparsity of~$\bfbeta$ implies the sparsity of $\bftheta$.
We assume that
\begin{equation}\label{eqc2}
\| \bftheta\| = O(n^{\tau} ), \qquad\tau< \eta
\mbox{ and }\tau\mbox{ does not depend on }n.
\end{equation}

If the number of nonzero components of~$\bfbeta$ is
$O(n^{2\tau})$, and all absolute values of
nonzero components of~$\bfbeta$ are bounded
by a constant $M$, then (\ref{eqc2}) holds since $\| \bftheta\| \leq
\| \bfbeta\| \leq M n^\tau$.

\begin{theorem}\label{th1}
Assume model (\ref{model}) and
conditions (\ref{eqc1}) and (\ref{eqc2}).

\begin{longlist}[(ii)]
 \item[(i)] As $n \rightarrow\infty$,
$ E ( \bfl' \hatbftheta- \bfl' \bftheta)^2 = O( h_n^{-1} ) +
O(h_n^2 n^{-2(\eta- \tau)}) $ uniformly over
$p$-dimensional deterministic vector
$\bfl$ with $\| \bfl\| =1$.

\item[(ii)] $n^{-1} E \| \bfX\hatbftheta- \bfX\bftheta\|^2 = O(r_n n^{-1} )
+ O(h_n^2 n^{-(1+\eta- 2 \tau)})$.
\end{longlist}
\end{theorem}

Note that these results hold without any condition on the dimension $p$.
Theorem~\ref{th1}(i) shows that the mean squared error of $\bfl' \hatbftheta$
converges to 0 uniformly in $\bfl$ if $h_n \rightarrow\infty$ and $h_n
n^{-(\eta- \tau)}
\rightarrow0$. Theorem~\ref{th1}(ii)
gives the convergence
rate of the expected
$L_2$-norm error $E\| \bfX\hatbftheta- \bfX\bftheta\|^2$
for estimating $E( \bfy) = \bfX\bftheta$.
Since the dimension of $\bfX\bftheta$
is $n$, we say that an estimator~$\hatbfvartheta$ of~$\bftheta$ is $L_2$-consistent if
$n^{-1} E\| \bfX\hatbfvartheta- \bfX\bftheta\|^2 \rightarrow0$
as $n \rightarrow\infty$. Typically, $r_n/n$ does not converge
to 0 and, hence, $\bfX\hatbftheta$ may not be $L_2$-consistent.

To elaborate the motivation of using
the expected $L_2$-norm error $E \| \bfX\hat{\bfvartheta}- \bfX\bftheta\|^2$
as a performance measure
for an estimator $ \hat{\bfvartheta}$ of $\bftheta$, we consider
the problem of predicting future $y$-values on deterministic $\bfX$.
Let $\bfy_*$ be independent of $\bfy$ but with
the same distribution as $\bfy$.
For deterministic $\bfX$, it is typical
to assess the accuracy of the prediction $ \bfX\hat{\bfvartheta}$
using the average prediction mean squared error
$ n^{-1} E \| \bfy_* - \bfX\hat{\bfvartheta}\|^2$. It turns out that
\[
n^{-1} E \| \bfy_* - \bfX\hat{\bfvartheta}\|^2
= \sigma^2 + n^{-1} E \| \bfX\hat{\bfvartheta}
- \bfX\bftheta\|^2 .
\]
Hence, having a small expected $L_2$-norm error
is equivalent to having a small average prediction
mean squared error.

\section{The thresholded ridge regression estimator}\label{sec4}

The discussion in the previous section
indicates that, although the ridge regression estimator
$\hat{\bftheta}$ is consistent for the estimation of any linear
combination of $\bftheta$, it may not be $L_2$-consistent, that is,
$n^{-1} E\| \bfX\hatbftheta- \bfX\bftheta\|^2 $
may not converge to 0. To achieve $L_2$-consistency
(and good prediction property) under some sparsity
conditions on $\bftheta$, we propose to improve
the ridge regression estimator by thresholding.

\subsection{Variable selection}\label{sec4.1}
Let $\calM_{\bfbeta, 0}$ be the set of indices
of nonzero components of~$\bfbeta$, and let
$\widehat\calM$ be the set of indices of components of~$\bfbeta$
selected using a variable selection method.
The variable selection method or $\widehat\calM$ is
said to be selection-consistent if and only if
%
%
\renewcommand{\theequation}{\arabic{equation}}
\setcounter{equation}{6}
\begin{equation}
\lim_{n \rightarrow\infty}
P( \widehat\calM= \calM_{\bfbeta, 0} ) = 1. \label{sc}
\end{equation}
Unlike the case of random $\bfX$, for deterministic $\bfX$ with
$p>n$, the selection-consistency defined by
(\ref{sc}) is generally not achievable because~$\bfbeta$ is not
identifiable.
Some selection-consistency results for the case of $p>n$ and
deterministic $\bfX$ published in the literature
are based on very strong and sometimes unrealistic
conditions on the design matrix $\bfX$ to ensure the
identifiability of~$\bfbeta$.
In fact, when~$\bfbeta$ is not identifiable, it is not appropriate
to use~$\bfbeta$ to describe usefulness of components of $\mathbf{x}_i$, since
two different~$\bfbeta$ may result in the same responses
under model (\ref{model}). Although components of $\mathbf{x}_i$ corresponding
to zero components of~$\bfbeta$
are not related to $y_i$, due to the fact that
$\bfbeta$ is unknown and not identifiable,
these components of $\mathbf{x}_i$ may still be useful in
statistical analysis since we have to use model (\ref{model1})
instead of model (\ref{model}), that is, $\bftheta$
instead of~$\bfbeta$.

The previous discussion leads to variable selection in terms of
the projection vector $\bftheta$, since
any linear combination
$\bfl' \bfbeta$ is estimable if and only if
$\bfl' \bfbeta= \bfl' \bftheta$.
However, when
$\bfbeta$ contains many zero components, $\bftheta$ may not
have any zero component, although
many components of $\bftheta$ may be close to zero.
Small but not exactly zero components of $\bftheta$
do not contribute much in estimation
but add variability. Thus, we would like to carry out variable selection
in a more general sense as defined by Zhang and Huang (\citeyear{ZhaHua08}), that is,
we try to eliminate small components of $\bftheta$.
Condition (\ref{eqc4}) stated later may be used to define
whether a component of $\bftheta$ can be treated as small.\looseness=-1

We propose to threshold
the ridge regression estimator $\hat{\bftheta}$.
Let $\hat\theta_j$ be the $j$th components
of $\hatbftheta$, $j=1,\ldots,p$.
The thresholded ridge regression estimator is defined as~$\tilde{\bftheta}$ whose
$j$th component $\tilde\theta_j = \hat\theta_j$ if $|\hat\theta_j
| > a_n$
and $\tilde\theta_j =0$ if $|\hat\theta_j | \leq a_n$, $j=1,\ldots
,p$, where
%
\begin{equation}
a_n = C_1 n^{-\alpha},\qquad 0< \alpha\leq1/2 ,  C_1>0 ,
\label{an}
\end{equation}
is the thresholding value with $\alpha$ and $C_1$ not depending on $n$.
The computation of~$\tilde{\bftheta}$ is easy since it has an explicit form.
Thresholding can be viewed as a variable selection procedure; that is,
we select components of $\bftheta$ with indices in
$\calM_{\hat{\bftheta}, a_n}$, the set of indices of nonzero components
of $\tilde{\bftheta}$.
We now study the asymptotic behavior of $\calM_{\hatbftheta, a_n}$
under some conditions\vspace*{1pt} and appropriate choices of $a_n$ and $h_n$.
A condition on the divergence rate of $p=p_n$ as $n \rightarrow\infty
$ is
\renewcommand{\theequation}{C\arabic{equation}}
\setcounter{equation}{2}
\begin{equation}\label{eqc3}
p = O(e^{n^\nu} ), \qquad 0< \nu<1
\mbox{ and }\nu\mbox{ does not depend on }n.
\end{equation}
If $p = e^{n^\nu}$, it is referred to as the ultra-high dimension
[Fan and Lv (\citeyear{FanLv10})].\vadjust{\goodbreak}

\begin{theorem}\label{th2}
Assume model (\ref{model}) with normally distributed
$\varepsilon_i$ and conditions~(\ref{eqc1})--(\ref{eqc3}). Let $a_n$ be given by
(\ref{an}) with $\alpha< (\eta- \nu- \tau)/3$, $u_n =
1+ (\log\log n )^{-1}$
and $h_n = C_2 a_n^{-2}(\log\log n )^3 \log(n\vee p) $,
where $C_2>0$ is a~constant and $n \vee p$ is the maximum of $n$ and $p$.
Then, for any constant $t >0$,
%
\renewcommand{\theequation}{\arabic{equation}}
\setcounter{equation}{8}
\begin{equation}
P( \calM_{\bftheta, a_nu_n} \subset
\calM_{\hat{\bftheta}, a_n} \subset\calM_{\bftheta, a_n/u_n} )
= 1- O\bigl((n \vee p)^{-t}\bigr), \label{result1}
\end{equation}
where $\calM_{\bfxi, c_n}$ denotes the set of indices of components of
$\bfxi$
whose absolute values are larger than $c_n$.
\end{theorem}

Result (\ref{result1}) shows that, by thresholding
$\hat{\bftheta}$, we can eliminate all components of $\bftheta$
with absolute values less than $a_n/u_n$, but
keep all components of $\bftheta$
with absolute values larger than $a_nu_n$, with
probability tending to 1 at
the rate of $O((n\vee p)^{-t})$ for any $t>0$. This
rate is at least $O(n^{-t}) $ for any $t>0$ and
it is $O(e^{ - t n^\nu})$ for any $t >0$
if $ \log p$ has exactly the order $n^\nu$.

Let $q_{n-}$ and $q_{n+}$ be
the numbers of elements in $\calM_{\bftheta, a_nu_n}$
and $\calM_{\bftheta, a_n/u_n}$, respectively.
Then $q_{n-} \leq q_{n+}$.
Since $u_n \rightarrow1$, it is often true that
$q_{n+} - q_{n-} \rightarrow0$ as $n \rightarrow\infty$.
Then, result (\ref{result1}) implies that
%
\renewcommand{\theequation}{\arabic{equation}}
\setcounter{equation}{9}
\begin{equation}
P (\calM_{\hat{\bftheta}, a_n}= \calM_{\bftheta, a_n} )
= 1- O\bigl( (n\vee p)^{-t}\bigr), \label{result1a}
\end{equation}
which will be referred to as the consistency of $ \calM_{\hat{\bftheta}
, a_n}$.
This consistency is weaker than the selection-consistency
given by (\ref{sc}), but the latter may not be
achieved.

We now consider nonnormal $\varepsilon_i$
under model (\ref{model}), that is, the normality
assumption on $\varepsilon_i$ is replaced by
\renewcommand{\theequation}{M}
\begin{equation}\label{eqM}
E(\varepsilon_i^{k}) < \infty \qquad
\mbox{for an even integer $k$ not depending on $n$},
\end{equation}
and condition (\ref{eqc3}) is replaced by
\renewcommand{\theequation}{C\arabic{equation}$^{\prime}$}
\setcounter{equation}{2}
\begin{equation}\label{eqcprime}
p = O(n^l ),\qquad  1 \leq \mbox{$l< k/6$ and
$l$ does not depend on $n$,}
\end{equation}
while the other conditions, (\ref{eqc1}) and (\ref{eqc2}), remain the same.
When the normality condition is relaxed to the moment condition (\ref{eqM}),
we cannot handle a dimension at the divergence rate given by (\ref{eqc3}),
although the
polynomial-type divergence rate given by (\ref{eqcprime}) can still be much
larger than $n$.
The integer $k$ in condition (\ref{eqM}) has to be
sufficiently large so that $3l(t+1)/k < \eta- \tau$,
where $t>0$ is in the convergence rate of ${\calM}_{\hat{\bftheta}, a_n}$.

\renewcommand{\thetheorem}{\arabic{theorem}A}
\setcounter{theorem}{1}
\begin{theorem}\label{thm2A}
Assume model (\ref{model})
and conditions (\ref{eqM}), (\ref{eqc1}), (\ref{eqc2})\break
and~(\ref{eqcprime}). For any $t>0$, let $a_n$ be given by
(\ref{an}) with $\alpha\leq(\eta- \xi- \tau)/3$
and $\xi= 3l(t+1)/k$, and
$u_n = 1+ (\log\log n )^{-1}$.
If $h_n = C_2 a_n^{-2}(\log\log n )^2 (n\vee p)^{2\xi/(3l)} $, where
$C_2>0$ is a constant,
then result (\ref{result1}) holds.
\end{theorem}

\subsection{$L_2$-consistency}\label{sec4.2}

The following result shows that, after the variable selection,
the thresholded estimator $\tilde{\bftheta}$ has asymptotically
smaller expected
$L_2$-norm error than $\hat{\bftheta}$, and it is in fact\vadjust{\goodbreak}
$L_2$-consistent, under the following sparsity condition on $\bftheta$:
\renewcommand{\theequation}{C\arabic{equation}}
\setcounter{equation}{3}
\begin{equation}\label{eqc4}
q_{n+}-q_{n-} \rightarrow0,\qquad   q_n/r_n \rightarrow0\quad \mbox{and}\quad
a_n v_n \rightarrow0,
\end{equation}
where
\[
v_n = \sum_{j\dvtx  |\theta_j| \leq a_n} | \theta_{j} |,
\]
$\theta_j$ is the $j$th component of $\bftheta$,
$r_n$ is the rank of $\bfX$,
$a_n$ is given by (\ref{an}), and~$q_n$, $q_{n-}$ and $q_{n+}$ are,
respectively, the numbers of
elements in sets $\calM_{\bftheta, a_n}$,
$\calM_{\bftheta, a_nu_n}$ and
$\calM_{\bftheta, a_n/u_n}$ given by (\ref{result1}).

The last two
conditions in (\ref{eqc4}) are very similar to condition (2.4) in
Zhang and Huang (\citeyear{ZhaHua08}); that is, there exist
$q_n $ ``large'' components of $\bftheta$ with $q_n$
much smaller than the rank of $\bfX$, and $v_n$,
the $L_1$ norm of the ``small'' components
of $\bftheta$, may diverges to $\infty$, but at
a rate slower than $a_n^{-1}$.

\renewcommand{\thetheorem}{\arabic{theorem}}
\setcounter{theorem}{2}
\begin{theorem}\label{th3}
Assume the conditions in Theorem~\ref{th2}
or~\ref{thm2A}. Assume further that~(\ref{eqc4}) holds and
the maximum eigenvalue of $\bfX' \bfX$ is $O(n)$. Then
%
\renewcommand{\theequation}{\arabic{equation}}
\setcounter{equation}{10}
\begin{equation}
n^{-1} E \| \bfX\tilde{\bftheta} - \bfX\bftheta\|^2 = O(q_n n^{-1})
+ O(v_na_n) + O\bigl( h_n^2 n^{-(1+\eta- 2\tau)}\bigr).\label{result2}
\end{equation}
\end{theorem}

Result (\ref{result2}) shows the gain of variable selection
by thresholding.
The expected $L_2$-norm error
$n^{-1} E\| \bfX\tilde{\bftheta}- \bfX\bftheta\|^2$
is smaller than $n^{-1} E\| \bfX\hat{\bftheta}- \bfX\bftheta\|^2$ for
sufficiently
large $n$. The former converges to 0 at a certain rate and hence~$\tilde{\bftheta}$ is $L_2$-consistent, whereas
the latter may not converge to 0 when $r_n /n$ does not converge to 0.

If $q_n/n \rightarrow0$,
result (\ref{result2}) can also be established
with the vector of nonzero components
of $\tilde{\bftheta}$ replaced by
the ordinary least squares estimator
of the sub-vector of $\bftheta$
indexed by the set $\calM_{\hat{\bftheta},a_n}$.

\subsection{Tuning parameters}\label{sec4.3}

To apply thresholding, we need to choose the constants $C_1$ in
the thresholding value $a_n$ given by (\ref{an}) and $C_2$ in the
regularization parameter
$h_n$ given in Theorem~\ref{th2} or~\ref{thm2A}. Similar to many other problems,
$C_1$ and $C_2$ can be viewed as tuning parameters, and there is no
optimal way to
find their values. Some discussions can be found, for example, in Fan
and Lv (\citeyear{FanLv08}).
It is possible to use a data-driven method to find values of tuning
parameters by minimizing
the average prediction mean squared error $ n^{-1} E \| \bfy_* -
\bfX\tilde{\bftheta}\|^2 = \sigma^2 + n^{-1} E \| \bfX\tilde{\bftheta}
- \bfX\bftheta\|^2$.

Let $\psi(C)$ be the average prediction mean squared error
when $C= (C_1,C_2)$ is used in $a_n$ and $h_n$.
Since $\psi(C)$ is unknown, we minimize the cross-validation estimator
\[
\hat\psi(C) = \frac{1}{n} \sum_{i=1}^n \bigl(y_i -
\mathbf{x}_i' \tilde{\bftheta}_{-i}^{(C)}\bigr)^2 ,
\]
where $\tilde{\bftheta}_{-i}^{(C)}$ is
the thresholded ridge regression estimator
of $\bftheta$ based on the data set with
$(y_i,\mathbf{x}_i )$ removed, $i=1,\ldots,n$.
To avoid repeated computation of $\tilde{\bftheta}_{-i}^{(C)}$,
we may use an equivalent formula for $\hat\psi(C)$,
%
\renewcommand{\theequation}{\arabic{equation}}
\setcounter{equation}{11}
\begin{equation}
\hat\psi(C) = \frac{1}{n} \sum_{i=1}^n \biggl(\frac{y_i -
\mathbf{x}_i' \tilde{\bftheta}^{(C)}}{1-w_{i}(C)}\biggr)^2 , \label{cv}
\end{equation}
where $w_{i}(C) =
\mathbf{x}_i' (\bfX'\bfX+ h_nI_p)^{-1} \mathbf{x}_i$ and
$\tilde{\bftheta}^{(C)}$ is the thresholded
ridge regression estimator based on the whole data set.
This method is applied in the simulation study presented
in the next section.

\section{Simulation results}\label{sec5}

With deterministic $\bfX$ and $p>n$,
we examined the $L_2$-norm errors and the expected $L_2$-norm
errors of the ridge regression estimator, the
thresholded ridge regression estimator, and
the popular LASSO estimator and ENET estimator (for comparison purpose)
in four simulation studies. In the first two
simulation studies, the design matrix $\bfX$
was generated from a multivariate normal distribution
but fixed throughout the simulation, which
corresponds to analysis conditional on $\bfX$.
In the last two simulation studies, $\bfX$ is
a nearly orthogonal Latin hypercube design or a
Latin hypercube design.

\subsection{Simulation study I}\label{sec5.1}

\begin{figure}

\includegraphics{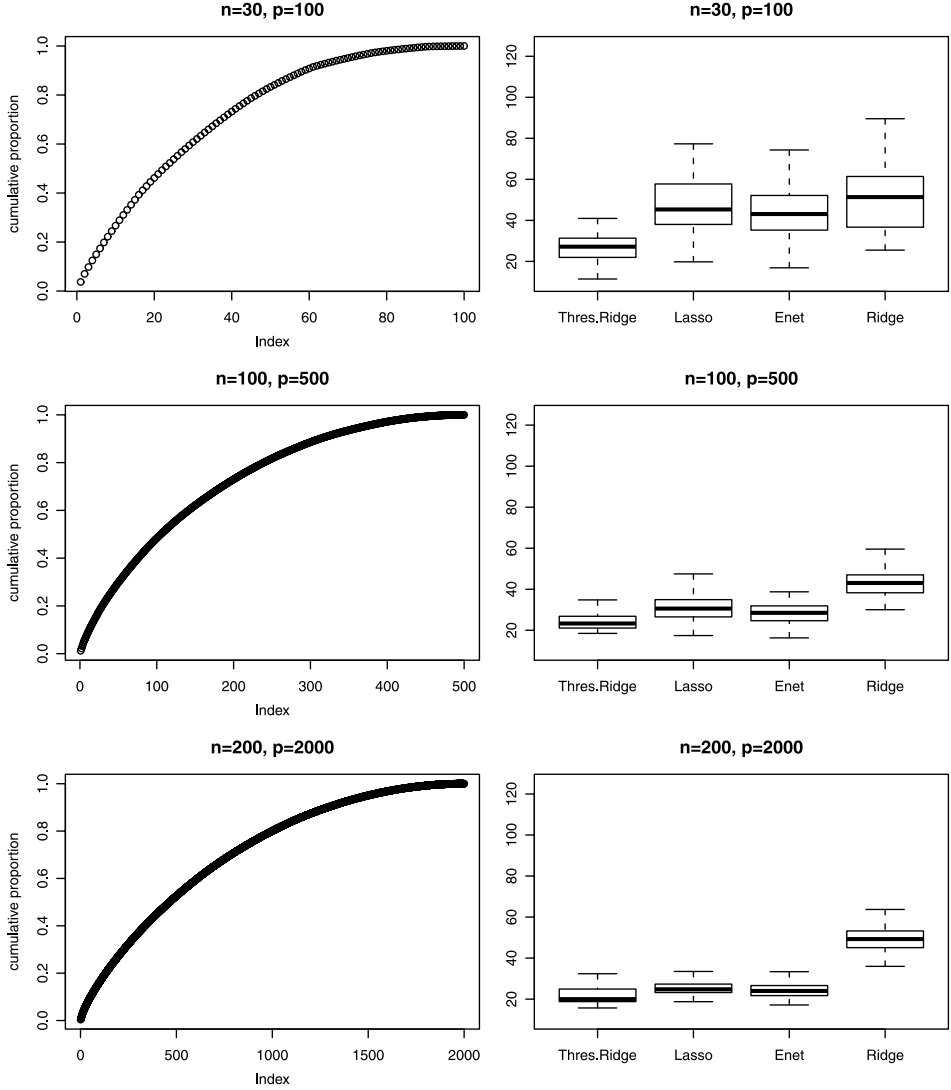}%
\vspace*{-3pt}
\caption{Study I: $L_2$ cumulative proportion plot
of $\bftheta$ and box plots of $L_2$-norm error for the thresholded ridge
regression, LASSO, ENET and ridge regression.}\label{fig1}
\vspace*{-3pt}
\end{figure}

We considered linear model (\ref{model}) with normally distributed
$\varepsilon_i $ and $\sigma= 10$.
Three sets of sample and variable sizes were considered,
$(n,p)= (30,100)$, $(100,500)$ and $(200,2000)$, with increasing ratio $p/n$.
A set of $\mathbf{x}_1,\ldots,\mathbf{x}_n$ were independently
generated with $\mathbf{x}_i \sim N( \bfzero, \bfSigma)$,
where the diagonal elements of $\bfSigma$ are all equal to 1
and off-diagonal elements of $\bfSigma$ are all equal to
0.75. This set of $\bfX$ was fixed throughout the simulation.
The first 20 components of~$\bfbeta$ are
$1+0.1j$ for $j=1,\ldots,20$, and the rest of the
components of~$\bfbeta$ are all equal to 0.
The $L_2$ cumulative proportion plot
of the projection vector $\bftheta$, that is,
$\sum_{j=1}^k \theta_{(j)}^2/\| \bftheta\|^2$, $k=1,\ldots,p$,
is given in Figure~\ref{fig1}, where $\theta_{(j)}^2$ is the
$j$th ordered value of $\theta_1^2,\ldots,\theta_p^2$. Although
$\bfbeta$
has many zero
components, $\bftheta$ does not have any zero component
but many components of $\bftheta$ are small.

%

For the thresholded ridge regression estimator,
we selected the tuning parameter $C=(C_1,C_2)$ by minimizing
$\hat\psi(C)$ given by (\ref{cv}). For the ridge
regression, LASSO, and ENET estimators, the tuning parameters
were selected by a 5-fold cross-validation.

Let $\hat{\bfvartheta}$ denote
the thresholded ridge regression estimator $\tilde{\bftheta}$,
the ridge regression estimator $\hat{\bftheta}$, the LASSO estimator or
the ENET estimator.
We independently generated 100 values of $\bfy$ and obtained 100
values of
$n^{-1} \| \bfX\bfbeta- \bfX\hat{\bfvartheta}\|^2$, the $L_2$-norm error
(divided by the sample size).
Box plots of 100 values of $n^{-1} \| \bfX\bfbeta- \bfX\hat{\bfvartheta}\|^2 $
for four estimation methods are given in Figure~\ref{fig1}.
The average of 100 values of $n^{-1} \| \bfX\bfbeta- \bfX\hat{\bfvartheta}\|^2 $,
a simulation approximation to the expected $L_2$-norm error
$n^{-1} E \| \bfX\bfbeta- \bfX\hat{\bfvartheta}\|^2 $,
is listed in Table~\ref{tab1} for each of the four methods.\vspace*{-3pt}

%
\begin{table}
\caption{Simulation approximation to the expected $L_2$-norm error}\label{tab1}
\vspace*{-3pt}
\begin{tabular*}{\textwidth}{@{\extracolsep{\fill}}ld{3.0}d{4.0}cccd{3.2}@{}}
\hline
& & & \multicolumn{4}{c@{}}{\textbf{Method}} \\[-6pt]
& & & \multicolumn{4}{c@{}}{\hrulefill} \\
\multicolumn{1}{@{}l}{\textbf{Study}} & \multicolumn{1}{c}{$\bolds{n}$} & \multicolumn{1}{c}{$\bolds{p}$} &
\multicolumn{1}{c}{\textbf{Thres. Ridge}} & \multicolumn{1}{c}{\textbf{LASSO}} & \multicolumn{1}{c}{\textbf{ENET}} & \multicolumn{1}{c@{}}{\textbf{Ridge}} \\
\hline
I & 30 & 100 & 27.34 & 48.46 & 44.56 & 51.48 \\
& 100 & 500 & 24.72 & 32.01 & 28.46 & 44.32 \\
& 200 & 2000 & 21.86 & 25.37& 24.17& 49.37\\[3pt]
II & 30 & 100 & 56.50 &69.05 &70.70 &76.05 \\
& 100 & 500 & 59.35 &68.33 &64.43 &94.06 \\
& 200 & 2000 & 74.59 &85.14 &82.35 &100.75 \\[3pt]
III & 49 & 96 & 61.58 &78.40 &76.83 &85.46 \\
& 64 & 192 & 54.79 &81.54 &79.78 &78.34 \\[3pt]
IV & 30 & 100 &43.44 &55.35 &49.29 &59.72\\
& 100 & 500 & 46.49 &56.60 &52.83 &65.85\\
& 200 & 2000 & 48.53 &51.78 &56.26 &71.21\\
\hline
\end{tabular*}
\vspace*{-3pt}
\end{table}

\subsection{Simulation study II}\label{sec5.2}

\begin{figure}

\includegraphics{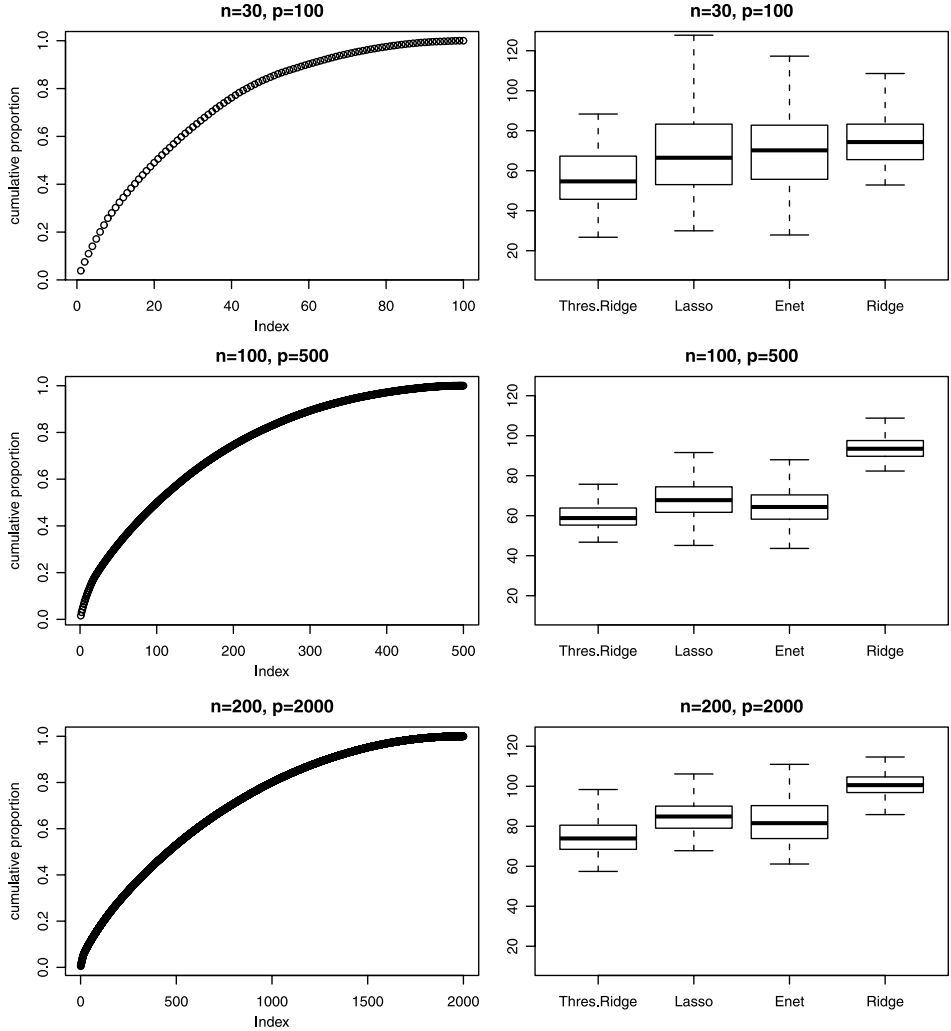}%
\caption{Study II: $L_2$ cumulative proportion plot
of $\bftheta$ and box plots of $L_2$-norm error for the thresholded ridge
regression, LASSO, ENET and ridge regression.}\label{fig2}
\end{figure}

The setting in this study is the same as that in
simulation study I except that the values
of $\mathbf{x}_i$'s were generated with
a~$\bfSigma$ whose
$(k,l)$th element is equal to
$(0.5)^{|k-l|}$ when $|k-l| \leq10$ and 0
when $|k-l| > 10$.\vadjust{\goodbreak}
The $L_2$ cumulative proportion plot
of $\bftheta$ and box plots of values of
$n^{-1} \| \bfX\bfbeta- \bfX\hat{\bfvartheta}\|^2 $
based on 100 simulation runs for four
estimation methods are given in Figure~\ref{fig2}.
The simulation approximations to $n^{-1} E \| \bfX\bfbeta- \bfX\hat{\bfvartheta}\|^2 $
are included in Table~\ref{tab1}.\vspace*{-3pt}

%

\subsection{Simulation study III}\label{sec5.3}

\begin{figure}

\includegraphics{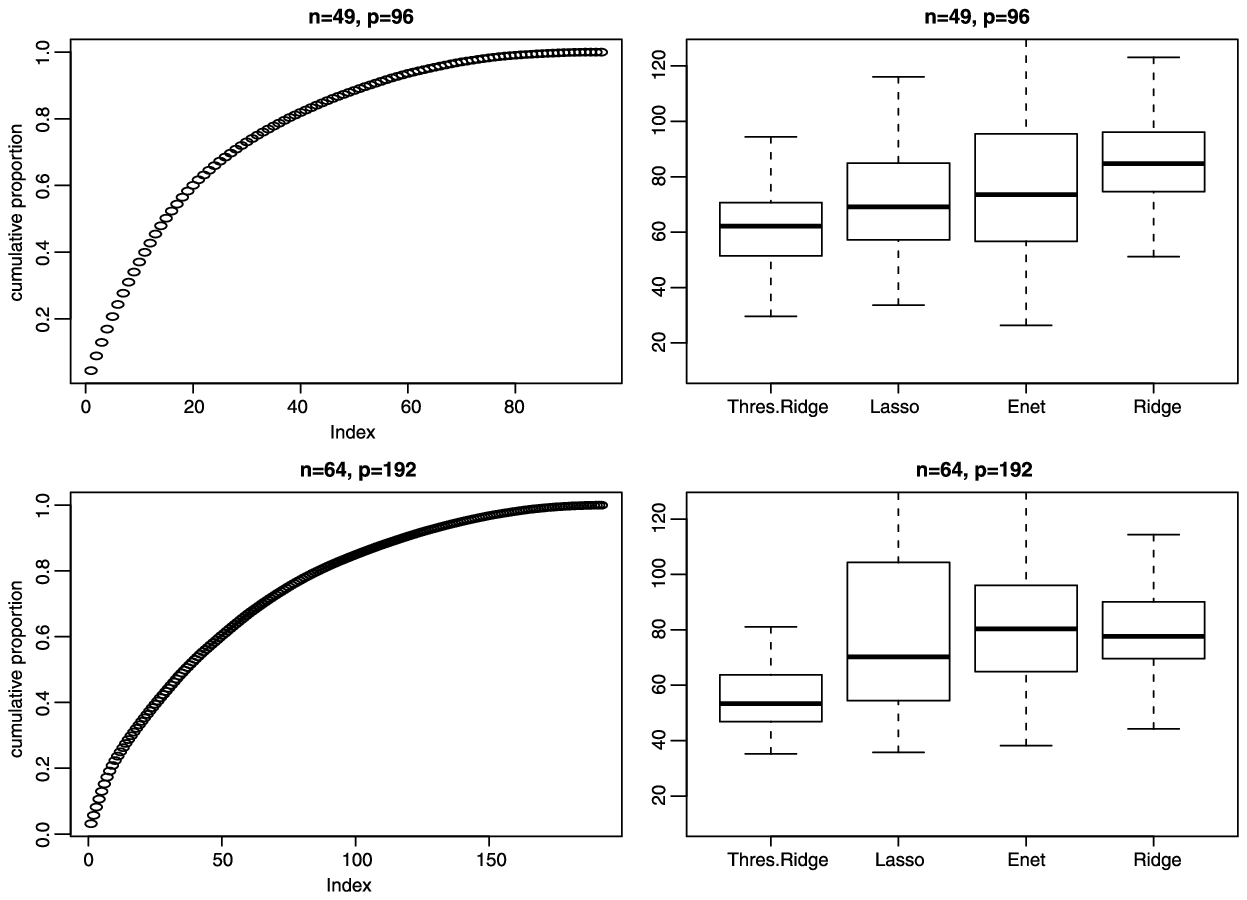}

\caption{Study III: $L_2$ cumulative proportion plot
of $\bftheta$ and box plots of $L_2$-norm error for the thresholded ridge
regression, LASSO, ENET and ridge regression.}\label{fig3}
\end{figure}

Let $\operatorname{NOLH}(n,p)$ denote a nearly orthogonal Latin
hypercube design with $n$ rows (runs) and $p$ columns
(variables). We considered two sets of $n$ and $p$.
In the first case, $n=49$, $p=96$ and
$\bfX$ is an
$\operatorname{NOLH}(49,96)$ constructed by using the orthogonal array-based
method in Lin, Mukerjee and Tang (\citeyear{LinMukTan09}).
In the second case, $n=64$, $p=192$ and
$\bfX$ is an $\operatorname{NOLH}(64,192)$. In both cases,
the first 15 components of~$\bfbeta$ are equal to
$0.2, 0.4,\ldots,2.8, 3.0$, and the rest components
of~$\bfbeta$ are equal to 0.
The standard deviation of $\varepsilon_i$ is $8$.
The rest of the simulation
setting is the same as that in simulation study I.
The $L_2$ cumulative proportion plot
of $\bftheta$ and box plots of values of
$n^{-1} \| \bfX\bfbeta- \bfX\hat{\bfvartheta}\|^2 $
based on 100 simulation runs for four
estimation methods are given in Figure~\ref{fig3}.
The simulation approximations to $n^{-1} E \| \bfX\bfbeta- \bfX\hat{\bfvartheta}\|^2 $
are included in Table~\ref{tab1}.\vspace*{-3pt}

%

\subsection{Simulation study IV}\label{sec5.4}

\begin{figure}

\includegraphics{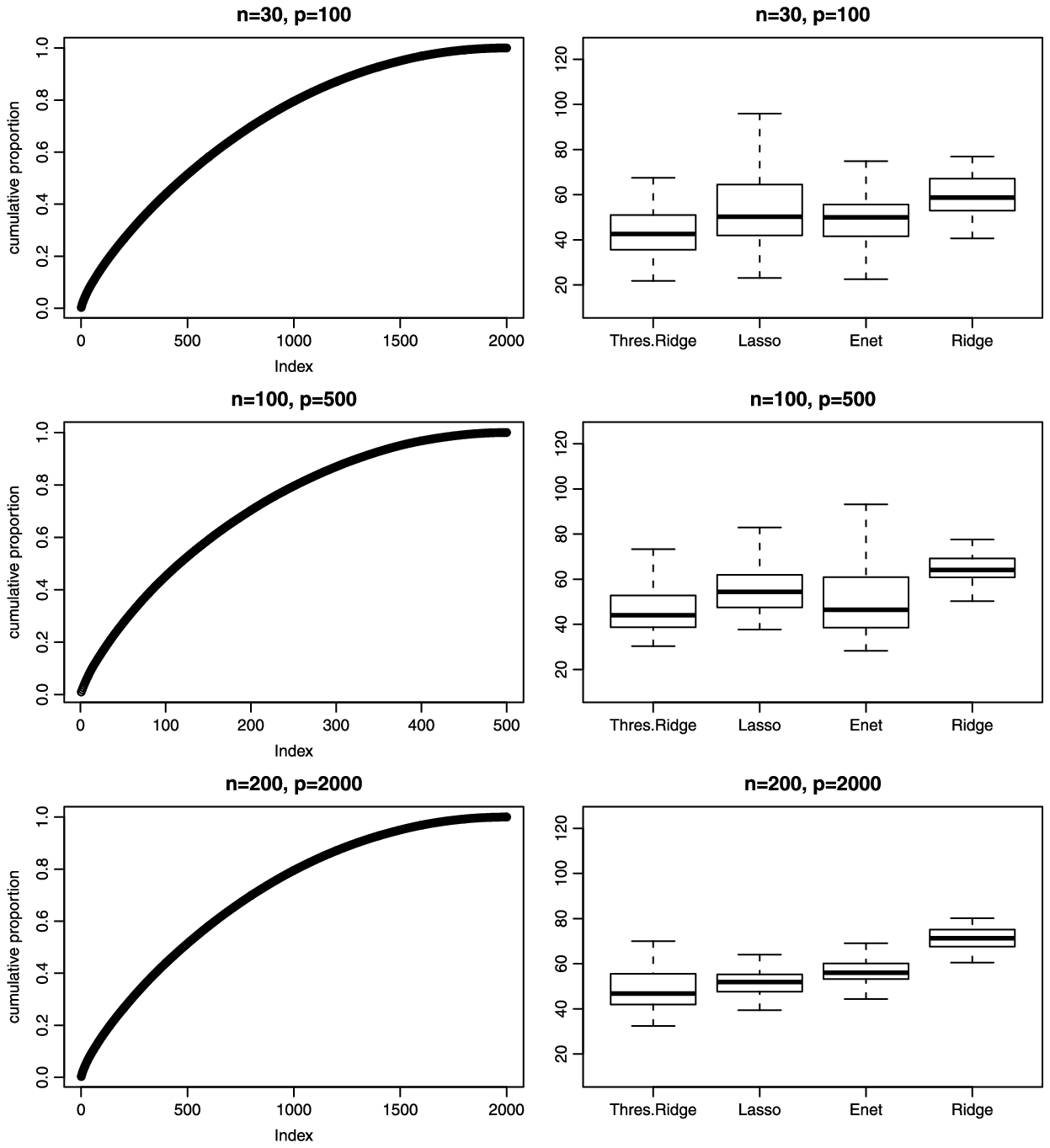}
\vspace*{-3pt}
\caption{Study IV: $L_2$ cumulative proportion plot
of $\bftheta$ and box plots of $L_2$-norm error for the thresholded ridge
regression, LASSO, ENET and ridge regression.}\label{fig4}
\vspace*{-3pt}
\end{figure}

The setting in this study is the same as that in
simulation study I except that $\bfX$ is a deterministic
Latin hypercube design [McKay, Beckman and Conover (\citeyear{McKBecCon79})]:
each column of $\bfX$ is a random permutation of~$n$ points $6(i/n) -3$, $i=1,\ldots,n$,
and all columns are generated independently.
The~$L_2$ cumulative proportion plot
of $\bftheta$ and box plots of values of
$n^{-1} \| \bfX\bfbeta- \bfX\hat{\bfvartheta}\|^2 $
based on 100 simulation runs for four
estimation methods are given in Figure~\ref{fig4}.
The simulation approximations to $n^{-1} E \| \bfX\bfbeta- \bfX\hat{\bfvartheta}\|^2 $
are included in Table~\ref{tab1}.\vspace*{-3pt}

%

\subsection{Conclusions based on simulation studies}\label{sec5.5}

From Table~\ref{tab1} and Figures~\ref{fig1}--\ref{fig4}, we conclude that the
thresholded ridge\vadjust{\goodbreak} regression estimator is much better than
the ridge regression estimator in terms of
the $L_2$-norm error or the expected $L_2$-norm error,
which supports our asymptotic theory, that is, the thresholded
ridge regression estimator is $L_2$-consistent whereas the ridge regression
estimator is not. Because the expected $L_2$-norm error is linearly
related to the average prediction mean squared error (Section~\ref{sec3}),
these results show that thresholding ridge regression has
better prediction performance. Except for study III, the LASSO
performs worse than the ENET and thresholded ridge regression, but
better than the ridge regression, and the ENET
performs worse than the thresholded ridge regression, although the difference
is small in some cases. Since the ENET uses a combination of $L_1$- and
$L_2$-penalty,
it is not surprising that its performance is between the LASSO and
thresholded ridge
regression. However, both LASSO and ENET have large variability
in simulation study III.
It is well known that the LASSO requires
more stringent conditions on the design matrix $\bfX$
[e.g., Zhao and Yu (\citeyear{ZhaYu06})].
The nearly orthogonal Latin hypercube design in simulation study III
may not satisfy these conditions, which results in the poor performance of
the LASSO. This also applies to the ENET, since it uses $L_1$-penalty.
Furthermore, no result for the $L_2$-consistency of LASSO or ENET
is available in the situation of deterministic $\bfX$ and $p>n$.

In terms of the computation, the thresholded ridge regression is much
simpler than
the LASSO or ENET, especially when $p$ is very large. Because of the identity~(\ref{iden}),
the computation complexity of the thersholded ridge regression
estimator does not increase as $p$ increases.

\section{Proofs}\label{sec6}
\mbox{}
\begin{pf*}{Proof of Lemma \protect\ref{le1}}
Suppose that (\ref{bspace}) holds. Let $\bfbeta_j \in\bfB$, $j=1,2$.
Then there are
$\bfxi_j \in\mathcal{R}^r$ such that $\bfbeta_j = \bfQ\bfxi_j +
\bfQ_\perp\phi( \bfxi_j )$, $j=1,2$. If $\bfX\bfbeta_1 = \bfX
\bfbeta_2$,
then, by (\ref{decom}), $\bfP\bfD\bfxi_1 = \bfP\bfD\bfxi_2$ and, thus,
$\bfxi_1 = \bfxi_2$, which implies $\bfbeta_1 = \bfbeta_2$.
This shows that the parameter~$\bfbeta$ in (\ref{model}) is identifiable.

Suppose now that $\bfB$ is not of the form (\ref{bspace}). Then,
there exist
$\bfxi\in\mathcal{R}^r$, $\bfzeta_j \in\mathcal{R}^{p-r}$, $j=1,2$,
$\bfzeta_1 \neq\bfzeta_2$ and $\bfbeta_j = \bfQ\bfxi+ \bfQ_\perp
\bfzeta_j \in\bfB$. Then $\bfbeta_1 \neq\bfbeta_2$,
but $\bfX\bfbeta_1 = \bfP\bfD\bfxi= \bfX\bfbeta_2$.
This shows that~$\bfbeta$ in (\ref{model}) is not identifiable.\vspace*{-3pt}
\end{pf*}

\begin{pf*}{Proof of Theorem \protect\ref{th1}}
\begin{longlist}[(ii)]
\item[(i)] From Section~\ref{sec3}, $\operatorname{bias}( \hatbftheta)
= - \bfQ
( h_n^{-1} \bfD^2 + \bfI_r )^{-1} \bfQ' \bftheta$.
From the facts that $ \bfQ'\bfQ=\bfI_r$,
$\bfD^2$ contains positive eigenvalues of $\bfX' \bfX$,\vadjust{\goodbreak}
and $( h_n^{-1} \bfD^2 + \bfI_r )^{-1}
\leq\frac{h_n/\lambda_{1n}}{1+h_n/\lambda_{1n}} \bfI_r$, we obtain that
$\| \operatorname{bias}( \hatbftheta)\| \leq\| \bftheta\| (h_n /
\lambda_{1n})$.
Hence, by (\ref{eqc1}) and (\ref{eqc2}),
$[ \bfl' \operatorname{bias} (\hatbftheta)]^2 \leq\| \operatorname
{bias} (\hatbftheta)\|^2$ $
= O( h_n^2 n^{-2(\eta- \tau)})$
uniformly over $\bfl$ with $\| \bfl\| =1$.
Also, from Section~\ref{sec3}, $\operatorname{var} ( \hatbftheta)
\leq\sigma^2 h_n^{-1} \bfI_p $.
Hence, $\bfl' \operatorname{var}(\hatbftheta) \bfl= O(h_n^{-1})$
uniformly over $\bfl$ with $\| \bfl\| =1$.
Then, the result follows from
$ E ( \bfl' \hatbftheta- \bfl' \bftheta)^2
= \bfl' \operatorname{var} ( \bftheta)\bfl+ [ \bfl' \operatorname
{bias} ( \bftheta)]^2 $.

\item[(ii)] Note that
$ E \| \bfX\hatbftheta- \bfX\bftheta\|^2 =
\operatorname{trace} [ \bfX\operatorname{var}( \hatbftheta) \bfX' ]
+ \| \bfX\operatorname{bias} ( \hatbftheta)\|^2 $.
From the proof of (i),
\begin{eqnarray*}
\bfX\operatorname{var}( \hatbftheta) \bfX' \back& \leq& \back
\sigma^2 \bfX( \bfX' \bfX+ h_n\bfI_p )^{-1} \bfX' \\
& = & \back\sigma^2 \bfP\bfD( \bfD^2 + h_n\bfI_r)^{-1} \bfD\bfP
' \\
& \leq& \back\sigma^2 \bfP\bfP' ,
\end{eqnarray*}
since $\bfD( \bfD^2 + h_n\bfI_r)^{-1} \bfD$ is a diagonal matrix whose
diagonal elements are bounded by 1.
Hence, $\operatorname{trace} [ \bfX\operatorname{var}( \hatbftheta)
\bfX' ]
\leq\sigma^2 \operatorname{trace}(\bfP\bfP' )$ $ = \sigma^2
r_n$. Also,
\[
\| \bfX\operatorname{bias} (\hatbftheta) \|^2 =
\bftheta' \bfQ( h_n^{-1} \bfD^2 + \bfI_r )^{-1}
\bfD^2 ( h_n^{-1} \bfD^2 + \bfI_r )^{-1} \bfQ' \bftheta
\leq h_n^2 \lambda_{1n}^{-1} \| \bftheta\|^2 ,
\]
which is $O( h_n^2 n^{-(\eta- 2 \tau)})$ by (\ref{eqc1}) and (\ref{eqc2}).
This completes the proof.\qquad\qed
\end{longlist}
\noqed\end{pf*}

\begin{pf*}{Proof of Theorem \protect\ref{th2}}
From the proof of Theorem~\ref{th1},
\[
\operatorname{bias}(\hat\theta_j ) = O ( \| \bftheta\| h_n /
\lambda_{1n} )
= O( h_n/ n^{\eta- \tau} )
\]
uniformly in $j=1,\ldots,p$.
For sufficiently large $n$, $\log\log n >0$.
With $h_n = C_2a_n^{-2}(\log\log n)^3 \log(n\vee p)$
and condition (\ref{eqc3}),
\[
\frac{h_n}{ n^{\eta- \tau} (u_n-1)a_n}
= \frac{ C_2(\log\log n)^4 \log(n\vee p)}{n^{\eta- \tau} a_n^3}
\leq\frac{c_1(\log\log n)^4}{n^{\eta- \nu- \tau- 3 \alpha}}
\]
for some constant $c_1>0$ and, hence,
$ |\operatorname{bias}(\hat\theta_j ) |/[(u_n-1)a_n] \rightarrow0 $
uniformly in $j$ when $\alpha<(\eta- \nu- \tau)/3$.
Since $\operatorname{var}(\hat\theta_j) = O(h_n^{-1})$,
there is a constant $c_0>0$ such that
\[
\frac{|\operatorname{bias}(\hat\theta_j ) |-
(u_n-1)a_n}{[\operatorname{var}(\hat\theta
_j)]^{1/2} }
\leq- \sqrt{2} c_0 \sqrt{h_n}a_n/(\log\log n ).
\]
Let $\Phi$ be the standard normal distribution function. From (\ref
{model}) with normally distributed $\varepsilon_i$,
\begin{eqnarray*}
P\bigl( | \hat\theta_j - \theta_{j} | > (u_n-1)a_n \bigr)
\back& \leq& \back
2 \Phi\biggl( \frac{|\operatorname{bias}(\hat\theta_j )| - (u_n-1)a_n }{
[\operatorname{var}(\hat\theta_j)]^{1/2}}\biggr)\\
& \leq& \back
2 \Phi\bigl( - \sqrt{2} c_0 \sqrt{h_n}a_n/(\log\log n) \bigr) \\
& \leq& \back\exp\{-c_0^2 h_n a_n^2/(\log\log n)^2 \} ,
\end{eqnarray*}
for sufficiently large $n$,
where the last inequality follows from $2 \Phi(-x) \leq e^{-x^2/2}$
for $x \geq2$
and the fact that $h_n a_n^2/(\log\log n)^2
= C_2 \log\log n \log(n\vee p) \rightarrow\infty$.
Using the same argument,
we also obtain that
\[
P\bigl( | \hat\theta_j - \theta_{j} | > (1-u_n^{-1})a_n \bigr)
\leq\exp\{-c_0^2 h_n a_n^2/(\log\log n)^2 \}
\]
for sufficiently large $n$.
Let $t>0$ be given. For sufficiently large $n$,
$c_0^2 C_2 \log\log n -1 >t$ and, hence,
\begin{eqnarray*}
P( \calM_{\bftheta, a_nu_n}
\subset{\calM}_{\hat{\bftheta}, a_n} ) \back& \geq& \back
1- P \biggl( \bigcup_{j \dvtx   | \theta_{j}| > u_na_n } \{ | \hat\theta_j |
\leq a_n \} \biggr) \\
& \geq& \back1-
P \biggl( \bigcup_{j \dvtx   | \theta_{j}| > u_na_n } \{ | \hat\theta_j - \theta_{j} | > (u_n-1)a_n
\} \biggr) \\
& \geq& \back
1 - \sum_{j=1}^p P \bigl( | \hat\theta_j - \theta_{j} | > (u_n-1)a_n
\bigr) \\
& \geq& \back
1- p \exp\{-c_0^2 h_n a_n^2/(\log\log n)^2 \} \\
& \geq& \back1 - (n\vee p)^{-t} .
\end{eqnarray*}
Similarly, for any $t>0$,
\begin{eqnarray*}
P( \calM_{\hat{\bftheta}, a_n} \subset\calM_{\bftheta, a_n/u_n})
\back& \geq& \back
P \biggl( \bigcap_{ j\dvtx   | \theta_j| \leq a_n /u_n} \{ | \hat\theta_j | \leq a_n \} \biggr) \\
& \geq& \back1-
P \biggl( \bigcup_{j \dvtx   | \theta_j| \leq a_n/u_n } \{ | \hat\theta_j - \theta_j | >
(1-u_n^{-1})a_n \} \biggr) \\
& \geq& \back
1- p \exp\{-c_0^2 h_n a_n^2/(\log\log n)^2 \} \\
& \geq& \back1 - (n\vee p)^{-t}
\end{eqnarray*}
for sufficiently large $n$. This completes the proof.
\end{pf*}

\begin{pf*}{Proof of Theorem \protect\ref{thm2A}}
From the proof of Theorem~\ref{th1}, we still have
$\operatorname{bias}(\hat\theta_j ) = O(h_n /n^{\eta-\tau})$ uniformly in
$j=1,\ldots,p$.
Let $\zeta_j$ be the $j$th component of
$(\bfX'\bfX+ h_n \bfI_p )^{-1} \sum_{i=1}^n \mathbf{x}_i
(y_i - \mathbf{x}_i' \bftheta)$. Then, for $u_n = 1+ (\log\log n)^{-1}$,
\begin{eqnarray*}
P\bigl( | \hat\theta_j - \theta_j | > (u_n-1)a_n \bigr)
\back& \leq& \back
\frac{E ( \hat\theta_j - \theta_j )^{k}}{[(u_n-1)a_n]^{k}} \\
& = & \back O\biggl( \frac{| \operatorname{bias}(\hat\theta_j )|^{k} +
E ( \zeta_j^{k})}{[(u_n-1)a_n]^{k}}\biggr) \\
& = & \back O\biggl( \frac{h_n^{k}(\log\log n)^{k}}{n^{k(\eta- \tau)}
a_n^{k}}\biggr)
+ O\biggl( \frac{(\log\log n)^{k}}{h_n^{k/2} a_n^{k}}\biggr) ,
\end{eqnarray*}
where the last equality follows from $E(\zeta_j^{k}) = O(h_n^{-k/2})$
[Whittle (\citeyear{Whi60}), Theorem~\ref{th2}].
Similarly,
\[
P\bigl( | \hat\theta_j - \theta_j | > (1-u_n^{-1})a_n \bigr)
= O\biggl( \frac{h_n^{k}(\log\log n)^{k}}{n^{k(\eta- \tau)}
a_n^{k}}\biggr)
+ O\biggl( \frac{(\log\log n)^{k}}{h_n^{k/2} a_n^{k}}\biggr).
\]
Using $h_n = C_2a_n^{-2}(\log\log n )^2 (n\vee p)^{2\xi/(3l)} $,
we obtain that
\begin{eqnarray*}
P( \calM_{\hat{\bftheta}, a_n} \subset\calM_{\bftheta, a_n/u_n}
)
\back& \geq& \back
1- \sum_{j=1}^p P\bigl( | \hat\theta_j - \theta_j | > (1-u_n^{-1})a_n
\bigr) \\[-2pt]
& = & \back
1 - O\biggl( \frac{p h_n^{k}(\log\log n)^{k}}{n^{k(\eta- \tau)}
a_n^{k}}\biggr)
- O\biggl( \frac{p (\log\log n)^{k}}{h_n^{k/2} a_n^{k}}\biggr)\\[-2pt]
& = & \back
1 - O\biggl( \frac{p(n\vee p)^{2k \xi/(3l)}(\log\log n)^{3k}}{n^{k(\eta
- \tau)}a_n^{3k}}\biggr)\\[-2pt]
& &{} \back-  O\biggl(\frac{p}{(n\vee p)^{ \xi k /(3l)}}\biggr)\\[-2pt]
& = & \back
1 - O\biggl( \frac{p(n\vee p)^{k \xi/l}(\log\log n)^{3k}}{
(n\vee p)^{(t+1)}n^{k(\eta- 3\alpha- \tau)}}\biggr) \\[-2pt]
& &{} \back-  O\biggl(\frac{p}{(n\vee p)^{(t+1)}}\biggr) \\[-2pt]
& = & \back
1 - O\biggl( \frac{n^{k \xi}(\log\log n)^{3k}}{(n\vee p)^{t}
n^{k(\eta- 3\alpha- \tau)}}\biggr) - O\biggl(\frac{1}{(n\vee
p)^{t}}\biggr)\\[-2pt]
& = & \back1 - o\bigl((n\vee p)^{-t}\bigr)- O\bigl((n\vee p)^{-t}\bigr) \\[-2pt]
& = & \back1 - O\bigl((n\vee p)^{-t}\bigr),
\end{eqnarray*}
since $k \xi/(3 l) =t+1$ and $\alpha\leq(\eta- \xi- \tau)/3$.
Similarly,
\[
P( \calM_{\bftheta, a_nu_n} \subset\calM_{\hat{\bftheta}, a_n}
) \geq1 -
O\bigl((n\vee p)^{-t}\bigr).
\]
Hence, result (\ref{result1}) follows.
\end{pf*}

\begin{pf*}{Proof of Theorem \protect\ref{th3}}
Let\vspace*{-1.5pt} $A_n = \{ \calM_{\hat{\bftheta},a_n} = \calM_{\bftheta, a_n} \}$
and $A_n^c$ be its complement. On the set $A_n$, the number of nonzero
components of $\tilde{\bftheta}$ is the same as $q_n$.
Let $\bftheta_1$ be $\bftheta$ with its components smaller than $a_n$ in
absolute value set to 0. Under condition (\ref{eqc4}) and the condition that
$\bfX' \bfX$ has a~maximum eigenvalue bounded by~$c n$ for a constant $c$,
\begin{eqnarray*}
n^{-1} \| \bfX\bftheta_1 - \bfX\bftheta\|^2 \back& \leq& \back
c \| \bftheta_1 - \bftheta\|^2 \\
& = & \back c \sum_{j\dvtx  | \theta_j | \leq a_n } \theta_j^2 \\
& \leq& \back c a_n \sum_{j\dvtx  | \theta_j | \leq a_n } |\theta_j| \\
& = & \back O(v_na_n ) .
\end{eqnarray*}
Hence,
\begin{eqnarray*}
n^{-1} E \| \bfX\tilde{\bftheta}- \bfX\bftheta\|^2
\back& \leq& \back2n^{-1} ( E\| \bfX\tilde{\bftheta}- \bfX
\bftheta_1 \|^2
+ \| \bfX\bftheta_1 - \bfX\bftheta\|^2 ) \\
& = & \back2 n^{-1} E \| \bfX\tilde{\bftheta}- \bfX\bftheta_1 \|^2
+ O(v_na_n) .
\end{eqnarray*}
Then, it remains to show that
%
\begin{equation}
n^{-1} E \| \bfX\tilde{\bftheta}- \bfX\bftheta_1 \|^2 =
O(q_n n^{-1})+ O(v_na_n) + O\bigl( h_n^2 n^{-(1+\eta- 2\tau)}\bigr). \label{result2a}
\end{equation}
Following the proof of Theorem~\ref{th1} we obtain that
\[
n^{-1} E[ \| \bfX\tilde{\bftheta}- \bfX\bftheta_1 \|^2 I_{A_n} ]
= O(q_n n^{-1})+ O\bigl( h_n^2 n^{-(1+\eta- 2\tau)}\bigr),
\]
where $I_A$ is the indicator of the set $A$. From
\[
\| \bfX\tilde{\bftheta}- \bfX\bftheta_1 \|^2 I_{A_n^c}
\leq2\| \bfX\tilde{\bftheta}- \bfX\hat{\bftheta}\|^2 I_{A_n^c}
+ 2\| \bfX\hat{\bftheta}- \bfX\bftheta_1 \|^2 I_{A_n^c}
\]
and Theorem~\ref{th1}, result (\ref{result2a}) follows if we can show that
\[
n^{-1} E\| \bfX\tilde{\bftheta}- \bfX\hat{\bftheta}\|^2 I_{A_n^c}
= o\bigl( q_nn^{-1} \vee h_n^2 n^{-(1+\eta- 2\tau)}\bigr).
\]
Since
\[
\| \bfX\tilde{\bftheta}- \bfX\hat{\bftheta}\|^2
= (\tilde{\bftheta}- \hat{\bftheta})' \bfX' \bfX( \tilde{\bftheta}-
\hat{\bftheta})
\leq O(n ) \| \tilde{\bftheta}- \hat{\bftheta}\|^2 \leq
O( a_n^2 p n ),
\]
the result follows from
$P( A_n^c ) = O((n\vee p)^{-t}) $ for any $t>0$ according to
Theorem~\ref{th2} or~\ref{thm2A}. This completes the proof.\qed
\noqed\end{pf*}

\section*{Acknowledgments}
The authors would like to thank a referee and an Associate Editor for
their helpful comments and suggestions.



\printaddresses

\end{document}